\documentclass[12pt]{amsart}
\usepackage{amscd,verbatim}
\usepackage[all]{xy}

\newcommand{\by}[1]{\overset{#1}{\longrightarrow}}
\newcommand{\iso}{\by{\sim}}

\newcommand{\C}{\underline{C}}

\newcommand{\Z}{\mathbf{Z}}

\newtheorem{thm}{Theorem}[section]
\newtheorem{lemma}[thm]{Lemma}
\newtheorem{prop}[thm]{Proposition}
\theoremstyle{definition}
\newtheorem{defn}[thm]{Definition}
\theoremstyle{remark}
\newtheorem{rk}[thm]{Remark}

\begin{document} 

\title{A note on relative duality for Voevodsky motives}
\author{Luca Barbieri-Viale}
\address{Dipartimento di Matematica Pura e Applicata, Universit\`a degli Studi di Padova\\ Via
Trieste, 63\\ I-35121 --- Padova\\ Italy}
\email{barbieri@math.unipd.it}
\author{Bruno Kahn}
\address{Institut de Math{\'e}matiques de Jussieu\\
175--179 rue du Chevaleret\\75013 Paris\\ France}
\email{kahn@math.jussieu.fr}
\date{January 23, 2007}
\maketitle

\section*{Introduction}

Relative duality is a useful tool in algebraic geometry and has been used several times. Here
we prove a version of it in Voevodsky's triangulated category of geometric motives $DM_{gm}(k)$
\cite{voetri}, where $k$ is a field which admits resolution of singularities. 

Namely, let $X$ be a smooth proper $k$-variety of pure dimension $n$ and $Y,Z$ two disjoint
closed subsets of $X$. We prove in Theorem \ref{t1} an isomorphism
\[M(X-Z,Y)\simeq M(X-Y,Z)^*(n)[2n]\]
where $M(X-Z,Y)$ and $M(X-Y,Z)$ are relative Voevodsky motives, see Definition \ref{d1}.

This isomorphism remains true after application of any $\otimes$-functor from $DM_{gm}(k)$,
for example one of the realisation functors appearing in \cite[I.VI.2.5.5 and I.V.2]{levine},
\cite{huber} or \cite{ivorra1}. In particular, taking the Hodge realisation, this
makes the recourse to M. Saito's theory of mixed Hodge modules unnecessary in
\cite[Proof of 2.4.2]{barsri}.

The main tools in the proof of Theorem \ref{t1} are a good theory of extended Gysin morphisms, readily deduced
from D\'eglise's work (Section \ref{1}) and Voevodsky's localisation theorem for motives with
compact supports \cite[4.1.5]{voetri}. This may be used for an alternative presentation of some
of the duality results of \cite[\S 4.3]{voetri} (see Remark \ref{r4}). The arguments seem
axiomatic enough to be transposable to other contexts.

We assume familiarity with Voevodsky's paper \cite{voetri}, and use its notation  throughout.

\section{Relative motives and motives with supports}

\begin{defn} \label{d1} Let $X\in Sch/k$ and $Y\subseteq X$, closed. We set
\begin{align*}
M(X,Y) &= \C_*(L(X)/L(Y))\\
M^Y(X)&= \C_*(L(X)/L(X-Y)).
\end{align*}
\end{defn}

\begin{rk} This convention is different from the one of D\'eglise in
\cite{deglise1,deglise2,deglise3} where what we denote by $M^Y(X)$ is written $M(X,Y)$ (and
occasionally $M_Y(X)$ as well).
\end{rk}

Note that $L(Y)\to L(X)$ and $L(X-Y)\to L(X)$ are monomorphisms, so that we
have
functorial exact triangles
\begin{gather}
\notag M(Y)\to M(X)\to M(X,Y)\by{+1}\\
\label{e3} M(X-Y)\to M(X)\to M^Y(X)\by{+1}.
\end{gather}

We can mix the two ideas: for $Y,Z\subseteq X$ closed, define
\[M^Z(X,Y) =\C_*(L(X)/L(Y)+L(X-Z)).\]

\begin{lemma}\label{l1} If $Y\cap Z=\emptyset$, the obvious map $M^Z(X)\to M^Z(X,Y)$ is
an
isomorphism, and we have an exact triangle
\[M(X-Z,Y)\to M(X,Y)\by{\delta} M^Z(X)\by{+1}.\qed
\]
\end{lemma}

\section{Extended Gysin}\label{1}

In the situation of Lemma \ref{l1}, assume that $Z$ is smooth of pure codimension $c$. F.
D\'eglise has then constructed a purity isomorphism
\begin{equation}\label{e2}
p_{Z\subset X}:M^Z(X)\iso M(Z)(c)[2c]
\end{equation}
with the following properties:

\begin{enumerate}
\item $p_{Z\subset X}$ coincides with Voevodsky's purity isomorphism of \cite[3.5.4]{voetri}
(see
\cite[1.11]{deglise3}).
\item If $f:X'\to X$ is transverse to $Z$ in the sense that $Z'=Z\times_X X'$ is smooth of pure
codimension $c$ in $X'$, then the diagram
\[\begin{CD}
M^{Z'}(X')@>p_{Z'\subset X'}>> M(Z')(c)[2c]\\
@V(f,g)_*VV @Vg_*VV\\
M^Z(X)@>p_{Z\subset X}>> M(Z)(c)[2c]
\end{CD}\]
commutes, where $g=f_{|Z'}$ (\cite[Rem. 4]{deglise1} or \cite[2.4.5]{deglise2}).
\item If $i:T\subset Z$ is a closed subset, smooth of codimension $d$ in $X$, the diagram
\[\xymatrix{
M^Z(X)\ar[r]^{p_{Z\subset X}}\ar[dd]^{i^*}
& M(Z)(c)[2c]\ar[dr]^{\alpha}\\
&& M^T(Z)(c)[2c]\ar[dl]_{p_{T\subset Z}}\\
M^T(X)\ar[r]^{p_{T\subset X}} &M(T)(d)[2d]
}\]
commutes, where $\alpha$ is the twist/shift of the boundary map in the triangle corresponding
to \eqref{e3} \cite[proof of 2.3]{deglise3}.
\end{enumerate}

\begin{defn}\label{d2} We set:
\[g_{Z\subset X}^Y=p_{Z\subset X}\circ \delta\]
where $p_{Z\subset X}$ is as in \eqref{e2} and $\delta$ is the morphism appearing in Lemma
\ref{l1}.
\end{defn}

In view of the properties of $p_{Z\subset X}$, these extended Gysin morphisms have the following
properties:

\begin{prop}\label{p2} a) Let $f:X'\to X$ be a morphism of smooth schemes. Let $Z'=f^{-1}(Z)$
and $Y'=f^{-1}(Y)$. If $f$ is transverse to $Z$, the diagram
\[\begin{CD}
M(X',Y')@>{g_{Z'\subset X'}^{Y'}}>> M(Z')(c)[2c]\\
@V{f_*}VV @V{g_*}VV\\
M(X,Y)@>{g_{Z\subset X}^Y}>> M(Z)(c)[2c]
\end{CD}\]
commutes, with $g=f_{|Z}$.\\
b) Let $X\supset Z\supset Z'$ be a chain of smooth $k$-schemes of pure
codimensions, and let $d=codim_Z Z'$. Let $Y\subset X$ be closed, with $Y\cap Z=\emptyset$. Then
\[g_{Z'\subset X}^Y = g_{Z'\subset Z}(d)[2d]\circ g_{Z\subset X}^{Y}.\]
\end{prop}

\section{Relative duality}

In this section, $X$ is a smooth proper variety purely of dimension $n$ and
$Y,Z$ are two disjoint closed subsets of $X$.  Consider the diagonal embedding
of $X$ into $X\times X$: its intersection with $(X-Y)\times (X-Z)$ is closed
and isomorphic to $X-Y-Z$. The closed subset $(X-Y)\times Y\cup Z\times(X-Z)$
is disjoint from $X-Y-Z$; from Definition \ref{d2} we get a extended
Gysin map
\begin{multline*}
M((X-Y)\times (X-Z),(X-Y)\times Y\cup Z\times(X-Z))\\
\to M(X-Y-Z)(n)[2n].
\end{multline*}

Note that the left hand side is isomorphic to $M(X-Y,Z)\otimes M(X-Z,Y)$ by an
explicit computation from the definition of relative motives. Composing with
the projection $M(X-Y-Z)(n)[2n]\to \Z(n)[2n]$, we get a map
\[M(X-Y,Z)\otimes M(X-Z,Y)\to \Z(n)[2n]\]
hence a map
\begin{equation}\label{e1}
M(X-Z,Y)\by{\alpha_X^{Y,Z}} M(X-Y,Z)^*(n)[2n].
\end{equation}

\begin{thm}\label{t1} The map \eqref{e1} is an isomorphism.
\end{thm}

The proof is given in the next section.

\section{Proof of Theorem \protect\ref{t1}}
 
\begin{lemma} \label{l2} If $Y=Z=\emptyset$ and $X$ is projective, then \eqref{e1} is an
isomorphism.
\end{lemma}

\begin{proof} As pointed out in \cite[p. 221]{voetri}, $\alpha_X^{\emptyset,\emptyset}$
corresponds to the class of the diagonal; then Lemma \ref{l2} follows from the functor of
\cite[2.1.4]{voetri} from Chow motives to $DM_{gm}(k)$. (This avoids a recourse to
\cite[4.3.2 and 4.3.6]{voetri}.)
\end{proof}

The next step is when $Z$ is empty. For any $U\in Sch/k$, write $M^c(U):=\C_*(L^c(U))$
\cite[p. 224]{voetri}. Since $X$ is proper, by \cite[4.1.5]{voetri} there is a canonical
isomorphism
\[M(X,Y)\iso M^c(X-Y)\]
induced by the map of Nisenvich sheaves
\[L(X)/L(Y)\to L^c(X-Y).\]

Therefore, from $\alpha_X^{Y,\emptyset}$, we get a map
\[\beta_X^Y:M^c(X-Y)\to M(X-Y)^*(n)[2n].\]

\begin{lemma}\label{l3}
The map $\beta_X^Y$ only depends on $X-Y$.
\end{lemma}

\begin{proof} Let $U=X-Y$. If $X'$ is another smooth compactification of $U$, with $Y'=X'-U$,
we need to show that $\beta_X^Y=\beta_{X'}^{Y'}$. By resolution of singularities, $X$ and $X'$
may be dominated by a third smooth compactification; therefore, without loss of generality, we
may assume that the rational map $q:X'\to X$ is a morphism. The point is that, in the diagram
\[\xymatrix{
M(X',Y')\ar[dr]\ar[ddr]_{\simeq}\ar[drr]^{\alpha_{X'}^{Y',\emptyset}}\\
& M(X,Y)\ar[r]_{\alpha_{X}^{Y,\emptyset}}\ar[d]^\simeq& M(U)^*(n)[2n]\\
& M^c(U)
}\]
both triangles commute. For the left one it is obvious, and for the upper one this follows
from the naturality of the pairing \eqref{e1}. Indeed, the square
\[\begin{CD}
X'-Y'@>\Delta'>> (X'-Y')\times X'\\
@V{q'}VV @V{q'\times q}VV\\
X-Y@>\Delta>> (X-Y)\times X
\end{CD}\]
is clearly transverse, where $q'=q_{|X'-Y'}$ (an isomorphism) and $\Delta,\Delta'$ are the
diagonal embeddings; therefore we may apply Proposition \ref{p2} a).
\end{proof}

From now on, we write $\beta_{X-Y}$ for the map $\beta_X^Y$.

\begin{lemma}\label{l4} a) Let $U\in Sm/k$ of pure dimension $n$, $T\by{i} U$ closed, smooth
of pure dimension $m$ and $V=U-T\by{j} U$. Then the diagram
\[\begin{CD}
M^c(T)@>\beta_T>> M(T)^*(m)[2m]\\
@V{i_*}VV @VV{g_{T\subset U}^*(n)[2n]}V\\
M^c(U)@>\beta_U>> M(U)^*(n)[2n]\\
@V{j^*}VV @VV{j^*}V\\
M^c(V)@>\beta_V>> M(V)^*(n)[2n]
\end{CD}\]
commutes.\\
b) Suppose that $\beta_T$ is an isomorphism. Then $\beta_U$ is an isomorphism if and only if
$\beta_V$ is.
\end{lemma}

\begin{proof} a) The bottom square commutes by a trivial case of Proposition \ref{p2} a). For
the top square, the statement is equivalent to the commutation of the diagram
\[\xymatrix{
&M^c(T)\otimes M(T)(c)[2c]\ar[dr]\\
M^c(T)\otimes M(U)\ar[ru]^{1\otimes g_{T\subset U}}\ar[rd]_{i_*\otimes 1} && \Z(n)[2n]\\
&M^c(U)\otimes M(U)\ar[ur]}\]
with $c=n-m$.

Take a smooth compactification $X$ of $U$, and let $\bar T$ be a desingularisation of the closure of $T$ in $X$. Let $q:\bar T\to X$ be the corresponding morphism, $Y=X-U$  and $W=\bar T-T$: we have to show that the diagram
\[\xymatrix{
&M(\bar T,W)\otimes M(T)(c)[2c]\ar[dr]\\
M(\bar T,W)\otimes M(U)\ar[ru]^{1\otimes g_{T\subset U}}\ar[rd]_{q_*\otimes 1} && \Z(n)[2n]\\
&M(X,Y)\otimes M(U)\ar[ur]}\]
or equivalently
\[\xymatrix{
&M(\bar T\times T,W\times T)(c)[2c]\ar[dr]\\
M(\bar T\times U,W\times U)\ar[ru]^{ f\circ g_{\bar T\times T\subset\bar T\times U}^{W\times U}}\ar[rd]_{(q\times 1)_*} && \Z(n)[2n]\\
&M(X\times U,Y\times U)\ar[ur]}\]
commutes, where $f$ is the map $M(\bar T\times T)(c)[2c]\to M(\bar T\times T,W\times T)(c)[2c]$.
For this, it is enough to show that the diagram
\[\xymatrix{
&M(\bar T\times T,W\times T)(c)[2c]\ar[r]_{g_{T\subset \bar T\times T}^{W\times T} (c)[2c]}&M(T)(n)[2n]\ar[dd]_{i_*}\\
M(\bar T\times U,W\times U)\ar[ru]^{ f\circ g_{\bar T\times T\subset\bar T\times U}^{W\times U}}\ar[rd]_{(q\times 1)_*} \\
&M(X\times U,Y\times U)\ar[r]^{g_{U\subset X\times U}^{Y\times U}}&M(U)(n)[2n]
}\]
commutes. Since extended Gysin extends Gysin, Proposition \ref{p2} a) shows that this amounts to
the commutatvity of 
\[\begin{CD}
M(\bar T\times U,W\times U)@>g_{T\subset \bar T\times U}^{W\times U}>> M(T)(n)[2n]\\
@V(q\times 1)_*VV @Vi_*VV\\
M(X\times U,Y\times U)@>g_{U\subset X\times U}^{Y\times U}>> M(U)(n)[2n]
\end{CD}\]
which follows from the functoriality of the extended Gysin maps (Proposition \ref{p2} b)).

b) This follows immediately from a).
\end{proof}

\begin{prop}\label{p3} $\beta_U$ is an isomorphism for all smooth $U$.
\end{prop}

\begin{proof} We argue by induction on $n=\dim U$, the case $n=0$ being known by Lemma
\ref{l2}. In general, let $V$ be an open affine subset of $U$ and pick a smooth projective
compactification $X$ of $V$, with $Z=X-V$. Let $Z\supset Z_1\supset\dots\supset Z_r=\emptyset$,
where $Z_{i+1}$ is the singular locus of $Z_i$. Let also $T=U-V$ and define similarly $T\supset
T_1\subset\dots\supset T_s=\emptyset$  (all $Z_i$ and $T_j$ are taken with their reduced
structure). Let $V_i=X-Z_i$ and $U_j=U-T_j$. Then $V_i-V_{i-1}$ and $U_j-U_{j-1}$ are smooth
for all $i,j$. Thus $\beta_U$ is an isomorphism by Lemma \ref{l2} (case of $\beta_X$) and a
repeated application of Lemma \ref{l4} b).
\end{proof}

\begin{rk} We haven't tried to check whether $\beta_U$ is the inverse of the isomorphism
appearing in the proof of \cite[4.3.7]{voetri}: we leave this interesting question to the
interested reader.
\end{rk}

\begin{proof}[End of proof of Theorem \ref{t1}]
 By Lemma \ref{l1}, the triangle $M(Z)\to
M(X-Y)\to M(X-Y,Z)\by{+1}$\ and the duality pairings  induce a map of
triangles
$$
\begin{CD}
M(X- Y, Z)^*(n)[2n]@>{}>>M(X- Y)^*(n)[2n]@>{}>>  M(Z)^*(n)[2n]\\
  @A{\alpha_X^{Y,Z}}AA @A{\alpha_X^{Y,\emptyset}}AA @A{\Phi}AA  \\
M(X-Z, Y)@>{}>>M(X, Y)@>{}>> M^Z(X).
\end{CD}
$$

(The left square commutes by a trivial application of Proposition \ref{p2} a), and $\Phi$ is
some chosen completion of the commutative diagram by the appropriate axiom of triangulated
categories.)

Consider the following diagram (which is the previous diagram with $Y=\emptyset$):
$$
\begin{CD}
M(X, Z)^*(n)[2n]@>{}>>M(X)^*(n)[2n]@>{}>>  M(Z)^*(n)[2n]\\
  @A{\alpha_X^{\emptyset,Z}}AA @A{\alpha_X^{\emptyset,\emptyset}}AA @A{\Phi}AA  \\
M(X-Z)@>{}>>M(X)@>{}>> M^Z( X)
\end{CD}
$$

Note that $\alpha_X^{\emptyset,Z}$ is dual to $\alpha_X^{Z,\emptyset}$; therefore it is an
isomorphism by Lemma \ref{l3} and Proposition \ref{p3}. It follows that $\Phi$ is an
isomorphism. Coming back to the first diagram and using Lemma \ref{l3} and Proposition \ref{p3}
a second time, we get the theorem.
\end{proof}

\begin{rk}
It would be interesting to produce a canonical pairing
\[
\cap_{(X, Z)}:M^Z( X)\otimes M(Z)\to \Z(n)[2n]
\]
playing the r\^ole of $\Phi$ in the above proof, i.e., compatible with $\alpha_X^{Y,Z}$.
\end{rk}

\begin{rk}\label{r4} As explained in \cite[App. B]{huka}, resolution of singularities and the
existence of the $\otimes$-functor of \cite[2.1.4]{voetri} are sufficient to prove that the
category $DM_{gm}(k)$ is rigid. Therefore, to apply the above arguments, one need only know
that the motives of the form $M(X-Y,Z)$ belong to $DM_{gm}(k)$, which is a consequence of
\cite[4.1.4]{voetri}.
\end{rk}

\end{document}